\documentclass[12pt,a4paper,twoside]{amsart}
\usepackage{geometry}
\usepackage{txfonts,pxfonts,tikz}
\usepackage{tgschola}
\usepackage{newtxmath}
\usepackage[T1]{fontenc}

\usepackage{amsmath}
\usepackage{amsthm}

\usepackage{mathtools}
\usepackage{stmaryrd}

\usepackage{wrapfig}
\usepackage{subcaption}

\usepackage{enumitem}
\usepackage{verbatim}
\usepackage{makecell}
\usepackage{xcolor}

\usepackage{hyperref}
\hypersetup{colorlinks,linkcolor={blue!75!black},citecolor={blue!75!black},urlcolor={blue!75!black}}

\usepackage{tikz}

\newcommand{\norm}[1]{\|#1\|}

\newcommand{\R}{\mathbb{R}}

\mathtoolsset{showonlyrefs}

\theoremstyle{theorem}

\newtheorem{theorem}{Theorem}[section] %

\theoremstyle{definition}

\theoremstyle{remark}

\numberwithin{equation}{section}

\newcommand{\abs}[1]{\left\vert#1\right\vert}

\newcommand{\dint}{\,\mathrm{d}}

\title[Sharp decay estimates for oscillatory integral operators]{Some sharp $L^2 \to L^p$ decay estimates for $(2+1)$-dimensional degenerate oscillatory integral operators}
\author{Shaozhen Xu}

\address{School of Information Engineering, Nanjing Xiaozhuang University, Nanjing 211171, China}
\email{shaozhen@nju.edu.cn}
\begin{document}
	
	\begin{abstract}
		We investigate $(2+1)-$dimensional oscillatory integral operators characterized by polynomial phase functions. By employing Stein's complex interpolation, we derive sharp $L^2\to L^p$ decay estimates for these operators.
	\end{abstract}
	
	\maketitle
	
	\setcounter{tocdepth}{1}
	
	\tableofcontents
	\section{Introduction}
	 The study of the asymptotic behavior of scalar oscillatory integrals of the form
	 \begin{equation}
	 	I(\lambda)=\int_{\R}e^{i\lambda f(x)}\phi(x)\dint x,
	 \end{equation}
     as the parameter $\lambda$ tends to infinity, is a fundamental problem in classical analysis with far-reaching applications, including those in chemical physics \cite{CCY92}. For nondegenerate phase functions $f(x)$, the stationary phase method yields an asymptotic expansion of the integral $I(\lambda)$. In the case of generic real-analytic phases, Varchenko's seminal work in \cite{Var76} established a profound connection between the principal term of the asymptotic expansion and the geometry of the Newton polyhedron associated with the phase function. A key ingredient in Varchenko’s proof is the application of toric resolution of singularities to desingularize the phase function. 
     
	 Like Fourier transform, the operator-analogue of oscillatory integrals with the form
	 	\begin{equation}
	 	\mathcal{T}_\lambda f(x)=\int_{\R^{n_{Y}}} e^{i\lambda S(x,y)}\psi(x,y)f(y)\dint y,\quad x\in\R^{n_X},
	 \end{equation}
	is called $\boldsymbol{\mathit{(n_X+n_Y)}}-$\textbf{\textit{dimensinal}} oscillatory integral operators.  We aim to establish the sharp $L^q\to L^p$ decay estimates for $\mathcal{T}_\lambda$. In the case $n_X=n_Y=n$ and $S(x,y)$ satisfies nondegenerate condition in the sense that the Hessian of $S(x,y)$ is non-vanishing in the support of $\psi(x,y)$, H\"{o}rmander \cite{Hor73} obtained the sharp $L^2\to L^2$ decay estimates with decay rate $-n/2$. In the degenerate case, we hope to describe the decay estimates for $\mathcal{T}_\lambda$ via Newton polyhedron as Varchenko had done for $I(\lambda)$. However, little is known except in the $(1+1)-$dimensional case, readers may refer to \cite{PhoSte92,See93,PhoSte94,PhoSte97,PhoSte98,PhoSteStu01,Ryc01,Yan04,Gre05,Xia17,SXY19} for more details. The increased complexity of higher-dimensional phase functions presents significant challenges, and the development of new techniques is required to address these difficulties. It should be noted that in \cite{Gre04} \cite{Xia17} resolution algorithm (Newton-Puiseux algorithm) for planar curves are used to deal with $\mathcal{T}_\lambda$. It is extremely difficult to explore whether or not higher dimensional resolution algorithms are suitable to  operator settings. 
	Further background and progress in high dimensions can be found in the \emph{Introduction} of \cite{Xu23} and the references therein.
	
	Fourier restriction operators, on one hand, can be regarded as a specific type of oscillatory integral operators subject to geometric conditions on their phase functions. This prompts us to explore the application of techniques originally developed for Fourier restriction problems to address degenerate oscillatory integral operators. In the historical context of the Fourier restriction conjecture, Fefferman resolved the planar case \cite{Fef70} using fractional integration. In higher dimensions, the Stein-Tomas estimate, derived through Stein's complex interpolation method, stood as the most powerful result until Bourgain \cite{Bou91} introduced Kakeya-type estimates to extend Fourier restriction estimates beyond the limitations of Stein-Tomas. It is noteworthy that the use of Kakeya-type estimates relies on the locally constant property of the Fourier restriction operator but does not apply to degenerate phases. Exploring the locally constant properties of degenerate oscillatory integral operators is an interesting and challenging direction for further research.

	In this paper, we aim to establish Stein-Tomas type estimates for the following $(2+1)-$dimensional oscillatory integral operators
	\begin{equation}\label{Key-O.I.O}
			T_\lambda f(x,y)=\int_{\R} e^{i\lambda (x^mt^k+y^nt^l)} \psi(x,y,t) f(t) \dint t,\quad (k,l,m,n)\in(\mathbb{Z}^{+})^4.
	\end{equation}
	The special case $m=1,k=2,n=2,l=1$ was addressed in \cite{Xu23} to obtain sharp $L^4\to L^4$ decay estimates using the broad-narrow method, and in \cite{TX24} to achieve sharp $L^2\to L^6$ decay using the fractional integration method. 
	
	In this paper, we shall establish the following $L^2\to L^p$ decay estimates for the operators \eqref{Key-O.I.O} by employing Stein's complex interpolation.
	\begin{theorem}\label{Main-Thm} 
		For the operator \eqref{Key-O.I.O}, without loss of generality, we may assume that $k>l$, then 
		\begin{equation}\label{Main-Ineq}
		\norm{T_\lambda f}_{L^{2k+2}(\R^2)}\leq C\lambda^{-\frac{1}{2(k+1)}\left(\frac{1}{m}+\frac{1}{\max\{n,l\}}\right)}\norm{f}_{L^2(\R)}.
		\end{equation}
		In particular, if further $l\leq n$, then the decay estimates above is sharp.
	\end{theorem}

\section{Proof of Theorem \ref{Main-Thm}}

Actually, we follow the main structure of Stein's proof in \cite{Ste93} with minor modifications. First, the $TT^{*}$ method enables us to reduce the inequality \eqref{Main-Ineq} to prove that $T_\lambda T_\lambda^{*}$ maps $L^{p'}(\R^2)$ into $L^p(\R^2)$. Second, we insert the operator $T_\lambda T_\lambda^{*}$ into a family of analytic operators. Third, we build $L^2(\R^2)\to L^2(\R^2)$ based on $(1+1)-$dimensional results of Phong and Stein \cite{PhoSte94}, and $L^1(\R^2)\to L^\infty(\R^2)$ based on van der Corput lemma. Last, Stein's complex interpolation yields the concluded inequality.

\begin{proof}
	We rewrite the inequality \eqref{Main-Ineq} in a dual version
	\begin{equation}\label{Dul-Main-Ineq}
		\norm{T_\lambda^{*} g}_{L^2(\R)}\leq C_{\psi}\lambda^{-\frac{1}{2(k+1)}\left(\frac{1}{m}+\frac{1}{\max\{n,l\}}\right)}\norm{g}_{L^{p'}(\R^2)}.
	\end{equation}
	where the dual opeator $T_\lambda^{*}$ is given by
	\begin{equation}
		T_\lambda^{*}g(t)=\int_{\R^2} e^{-i\lambda \left(u^mt^k+v^nt^l\right)}\bar{\psi}(u,v,t)g(u,v)\dint u\dint v.
	\end{equation}
	Notice that
	\begin{align}
		&\int_\R\left|T_\lambda^{*}g(t)\right|^2\dint t\\
		&=\int_\R T_\lambda^{*}g(t)\overline{T_\lambda^{*}g(t)}\dint t\\
		&=\int_\R\left[\int_{\R^2} e^{-i\lambda \left(u^mt^k+v^nt^l\right)}\bar{\psi}(u,v,t)g(u,v)\dint u\dint v\right] \cdot \left[\int_{\R^2} e^{i\lambda \left(x^mt^k+y^nt^l\right)}\psi(x,y,t)\bar{g}(x,y)\dint x\dint y\right]\dint t\\
		&=\int_{\R^2}\int_{\R^2}\left[\int_\R e^{i\lambda\left[(x^m-u^m)t^k+\left(y^n-v^n\right)t^l\right]}\bar{\psi}(u,v,t)\psi(x,y,t)\dint t \right] g(u,v)\bar{g}(x,y)\dint u\dint v\dint x\dint y.
	\end{align}
	If we set
	\begin{equation}
		K(x,y,u,v)=\int_\R e^{i\lambda\left[(x^m-u^m)t^k+\left(y^n-v^n\right)t^l\right]}\bar{\psi}(u,v,t)\psi(x,y,t)\dint t,
	\end{equation}
	then 
	\begin{equation}
		\int_\R\left|T_\lambda^{*}g(t)\right|^2\dint t=\int_{\R^2}\int_{\R^2}K(u,v,x,y) g(u,v)\bar{g}(x,y)\dint u\dint v\dint x\dint y.
	\end{equation}
	By writing 
	\begin{equation}
		T_Kg(x,y)=\int_{\R^2}K(u,v,x,y) g(u,v)\dint u\dint v,
	\end{equation}
	we have
	\begin{equation}\label{Dual-Red}
		\int_\R\left|T_\lambda^{*}g(t)\right|^2\dint t=\int_{\R^2}T_Kg(x,y) \bar{g}(x,y)\dint x\dint y.
	\end{equation}
	To prove \eqref{Main-Ineq}, it suffices to verify that
	\begin{equation}
		T_K:\quad L^{p'}\rightarrow L^{p}.
	\end{equation}
	To achieve this goal, we insert it into a family of analytic operators 
		\begin{equation}
		T_K^\alpha g(x,y)=\int_{\R^2}K^\alpha(u,v,x,y) g(u,v)\dint u\dint v,
	\end{equation}
	where the kernel
		\begin{equation}
		K^\alpha(u,v,x,y)=\int_{\R^2} e^{i\lambda\left[(x^m-u^m)t^k+\left(y^n-v^n\right)t^l+\phi(x,y,u,v,s)\right]}\bar{\psi}(u,v,t)\psi(x,y,t)\delta_\alpha(s)\dint t\dint s,
	\end{equation}
	satisfies
	\begin{equation}
		\phi(x,y,u,v,0)=0, \quad K^0(x,y,u,v)=K(x,y,u,v).
	\end{equation}
	
	Before presenting the formal construction of $\phi$ and $\delta_\alpha$, we provide some intuitive explanations. To build intuition, we may regard the function $\delta_\alpha$ as an analytic continuation of Dirac density $\delta$ in distributional sense and $\delta_0=\delta$, this implies
	\begin{equation}
		\left<\delta_0(s), f(s)\right>=\delta(f)=f(0), \quad \text{for any test function } f.
	\end{equation}
	
	While the construction of $\delta_\alpha$ is standard (see Chapter IX of \cite{Ste93}), we provide it here for completeness. We fix a $\zeta\in C_c^{\infty}$ with $\zeta(s)=1$ for $|s|\leq 1$, then
	\begin{equation}
		\delta_\alpha(s)=\begin{cases}
			\frac{e^{\alpha^2}}{\Gamma(\alpha)}s^{\alpha-1}\zeta^2(s), \quad &\text{if } s>0;\\
			
			0, &\text{if } s\leq 0.
		\end{cases}
	\end{equation}
	
	Next, we set
	\begin{equation}
		\phi(x,y,u,v,s)=(x^m-u^m)s.
	\end{equation}
	
	Write $K^\alpha$ explicitly,
	\begin{equation}
			K^\alpha(u,v,x,y)=\int_\R\int_\R e^{i\lambda\left[(x^m-u^m)(t^k+s)+\left(y^n-v^n\right)t^l\right]}\bar{\psi}(u,v,t)\psi(x,y,t)\delta_\alpha(s)\dint t\dint s.
	\end{equation}
	Then for $\mathrm{Re(\alpha)}=1$, 
	\begin{align}
		T_K^\alpha g(x,y)
		&=\int_{\R^2}K^\alpha(u,v,x,y) \bar{g}(u,v)\dint u\dint v\\
		&=\int_{\R^2}\int_\R\int_\R e^{i\lambda\left[(x^m-u^m)(t^k+s)+\left(y^n-v^n\right)t^l\right]}\bar{\psi}(u,v,t)\psi(x,y,t)\delta_\alpha(s)\dint t\dint s g(u,v)\dint u\dint v\\
		&:=S_1\circ S_2(g)(x,y),
	\end{align}
	where
	\begin{equation}
		S_1(h)(x,y)=\frac{e^{\alpha^2}}{\Gamma(\alpha)}\int_0^{+\infty}\int_\R e^{i\lambda\left[(x^m(t^k+s)+y^nt^l\right]}\psi(x,y,t)\zeta(s) s^{\alpha-1}h(t,s)\dint t\dint s,
	\end{equation}
	and 
		\begin{equation}
		S_2(g)(t,s)=\int_\R\int_\R e^{-i\lambda\left[(u^m(t^k+s)+v^nt^l\right]}\bar{\psi}(u,v,t)\zeta(s) g(u,v)\dint u\dint v.
	\end{equation}
	Since the coordinate transformation
	\begin{equation}
		w=s+t^k, t=t,
	\end{equation}
	is invertible everywhere, then iterating the $(1+1)-$dimensional result of Phong and Stein \cite{PhoSte94} implies 
	\begin{equation}\label{L2Est}
		\norm{T_K^\alpha g}_{L^2(\R^2)}\leq C\lambda^{-\frac{1}{m}-\frac{1}{\max\{n,l\}}}\norm{g}_{L^2(\R^2)} \quad Re(\alpha)=1.
	\end{equation}
	On the other hand, 
	\begin{equation}
		\norm{T_K^\alpha g}_{L^\infty(\R^2)}\leq C\sup\abs{K^\alpha(u,v,x,y)}\norm{g}_{L^1(\R^2)}.
	\end{equation}
	Observe that
	\begin{equation}
		K^\alpha(u,v,x,y)=\int_\R e^{i\lambda\left[(x^m-u^m)t^k+\left(y^n-v^n\right)t^l\right]}\bar{\psi}(u,v,t)\psi(x,y,t)\dint t\cdot \hat{\delta}_\alpha\left(\frac{\lambda(x^m-u^m)}{2\pi}\right).
	\end{equation}
	Van der Corput lemma tells 
	\begin{equation}
		\abs{\int_\R e^{i\lambda\left[(x^m-u^m)t^k+\left(y^n-v^n\right)t^l\right]}\bar{\psi}(u,v,t)\psi(x,y,t)\dint t}\leq C\left(1+\lambda\abs{x^m-u^m}\right)^{-1/k} \quad k>l.
	\end{equation}
	It remains to consider the upper bound of Fourier transform of $\delta_\alpha$,
	\begin{equation}
		\hat{\delta}_\alpha(t)=\frac{e^{\alpha^2}}{\Gamma(\alpha)}\int_{0}^{+\infty}e^{2\pi ist}s^{\alpha-1}\zeta^2(s)\dint s.
	\end{equation}
	To effectively evaluate oscillatory integrals of the form above, it is crucial to identify regions where the oscillatory phase function dominates the integrand's behavior and regions where the power function's characteristics prevail. Details can be found in \cite{Ste93}, we only state the conclusion
	\begin{equation}
		\abs{\hat{\delta}_\alpha(t)}\leq C\left(1+\abs{t}\right)^{-Re(\alpha)}.
	\end{equation}
	If we choose 
	\begin{equation}
		Re(\alpha)=-\frac{1}{k},
	\end{equation}
	then
	\begin{equation}
		\abs{K^\alpha(u,v,x,y)}\leq C,
	\end{equation}
	this implies
	\begin{equation}\label{L1Est}
			\norm{T_K^\alpha g}_{L^\infty(\R^2)}\leq C\norm{g}_{L^1(\R^2)}, \quad Re(\alpha)=-1/k.
	\end{equation}
	Combing \eqref{L2Est} with \eqref{L1Est}, Stein's complex interpolation yields
	\begin{equation}
		\norm{T_K g}_{L^{2k+2}(\R^2)}\leq C\lambda^{-\frac{1}{m(k+1)}-\frac{1}{\max\{n,l\}(k+1)}}\norm{g}_{L^{\frac{2k+2}{2k+1}}(\R^2)}.
	\end{equation}
	From \eqref{Dual-Red}, we know that
	\begin{equation}
		\norm{T_\lambda^{*} g}_{L^2(\R)}\leq C\lambda^{-\frac{1}{2(k+1)}\left(\frac{1}{m}+\frac{1}{\max\{n,l\}}\right)}\norm{g}_{L^{\frac{2k+2}{2k+1}}(\R^2)},
	\end{equation}
	which implies
	\begin{equation}
		\norm{T_\lambda f}_{L^{2k+2}(\R^2)}\leq C\lambda^{-\frac{1}{2(k+1)}\left(\frac{1}{m}+\frac{1}{\max\{n,l\}}\right)}\norm{f}_{L^2(\R)}.
	\end{equation}
	
To verify the optimality of the decay estimate for the case $n>l$, we refer to the example in \cite{Xu23} for demonstration.

Consider the function
	\begin{equation}\label{Supp-Func}
	\psi(x,y,t)=\begin{cases}
		0, &\quad \abs{(x,y,t)}\geq 1,\\
		1, &\quad \abs{(x,y,t)}\leq \frac{1}{2}.
	\end{cases}
\end{equation}
and choose
	\begin{equation}\label{Tes-Func}
		f(t)=\chi_{[0,1]}(t).
	\end{equation}
Assuming the prior decay estimate
	\begin{equation}
		\|T_
		\lambda f\|_{L^p(\R^2)}\leq C_\psi\lambda^{-\delta} \|f\|_{L^2(\R)},
	\end{equation}
	we obtain
	\begin{equation*}
		\left[\iint\abs{\int_0^1e^{i\lambda(x^mt^k+y^nt^l)}\psi(x,y,t)\dint t}^p\dint x\dint y\right]^{\frac{1}{p}}\leq C_\psi\lambda^{-\delta}.
	\end{equation*}
	On account of the support function \eqref{Supp-Func}, we have
	\begin{equation*}
		\lambda^{-\frac{1}{p}\left(\frac{1}{m}+\frac{1}{n}\right)}\lesssim \left[\int_{\abs{y}\lesssim \lambda^{-\frac{1}{n}}}\int_{\abs{x}\lesssim \lambda^{-\frac{1}{m}}}\abs{\int_0^1e^{i\lambda (x^mt^k+y^nt^l)}\psi(x,y,t)\dint t}^p\dint x\dint y\right]^{\frac{1}{p}}\leq C_\psi\lambda^{-\delta}.
	\end{equation*}
	Given that the inequality holds for arbitrarily large value of $\lambda$, it follows that
	\begin{equation*}
		\delta\leq \frac{1}{p}\left(\frac{1}{m}+\frac{1}{n}\right).
	\end{equation*}
	In the case where $p=2k+2$, this leads directly to our estimate, thereby completing the proof.
\end{proof}
	
\section*{Acknowledgments}
	The author is supported by Jiangsu Natural Science Foundation (Grant No. BK20200308 ), in part by the Research project of higher education institutions in Jiangsu Province (Grant No. 24KJB110021).


\begin{thebibliography}{SXY19}
		
		\bibitem[Bou91]{Bou91}
		J.~Bourgain.
		\newblock {$L^p$} estimates for oscillatory integrals in several variables.
		\newblock {\em Geom. Funct. Anal.}, 1(4):321--374, 1991.
		
		\bibitem[CY92]{CCY92}
		J.~L. Connor, P. R.~Curtis and W.~A. Young.
		\newblock Uniform asymptotics of oscillating integrals: applications in
		chemical physics.
		\newblock In {\em Wave Asymptotics, Proc. of the Meeting to Mark the Retirement
			of Fritz Ursell}, pages 283--310. Cambridge Univ. Press Cambridge, 1992.
		
		\bibitem[Fef70]{Fef70}
		C.~Fefferman.
		\newblock Inequalities for strongly singular convolution operators.
		\newblock {\em Acta Math.}, 124(1):9--36, 1970.
		
		\bibitem[Gre04]{Gre04}
		M.~Greenblatt.
		\newblock A direct resolution of singularities for functions of two variables
		with applications to analysis.
		\newblock {\em J. Anal. Math.}, 92:233--258, 2004.
		
		\bibitem[Gre05]{Gre05}
		M.~Greenblatt.
		\newblock Sharp {$L^2$} estimates for one-dimensional oscillatory integral
		operators with {$C^\infty$} phase.
		\newblock {\em Amer. J. Math.}, 127(3):659--695, 2005.
		
		\bibitem[H{\"o}r73]{Hor73}
		L.~H{\"o}rmander.
		\newblock Oscillatory integrals and multipliers on {F$L^p$}.
		\newblock {\em Ark. Mat.}, 11(1):1--11, 1973.
		
		\bibitem[PS92]{PhoSte92}
		D.~H. Phong and E.~M. Stein.
		\newblock Oscillatory integrals with polynomial phases.
		\newblock {\em Invent. Math.}, 110(1):39--62, 1992.
		
		\bibitem[PS94]{PhoSte94}
		D.~H. Phong and E.~M. Stein.
		\newblock Models of degenerate fourier integral operators and {Radon}
		transforms.
		\newblock {\em Ann. of Math.}, 140(3):703--722, 1994.
		
		\bibitem[PS97]{PhoSte97}
		D.~H. Phong and E.~M. Stein.
		\newblock The {Newton} polyhedron and oscillatory integral operators.
		\newblock {\em Acta Math.}, 179(1):105--152, 1997.
		
		\bibitem[PS98]{PhoSte98}
		D.~H. Phong and E.~M. Stein.
		\newblock Damped oscillatory integral operators with analytic phases.
		\newblock {\em Adv. Math.}, 134(1):146--177, 1998.
		
		\bibitem[PSS01]{PhoSteStu01}
		D.~H. Phong, E.~M. Stein, and J.~Sturm.
		\newblock Multilinear level set operators, oscillatory integral operators, and
		{Newton} polyhedra.
		\newblock {\em Math. Ann.}, 319(3):573--596, 2001.
		
		\bibitem[Ryc01]{Ryc01}
		V.~S. Rychkov.
		\newblock Sharp {$L^2$} bounds for oscillatory integral operators with
		{$C^\infty$} phases.
		\newblock {\em Math. Z.}, 236(3):461--489, 2001.
		
		\bibitem[See93]{See93}
		A.~Seeger.
		\newblock Degenerate fourier integral operators in the plane.
		\newblock {\em Duke Math. J.}, 71:685--745, 1993.
		
		\bibitem[SM93]{Ste93}
		E.~M. Stein and T.~S. Murphy.
		\newblock {\em Harmonic analysis: real-variable methods, orthogonality, and
			oscillatory integrals}, volume~3.
		\newblock Princeton University Press, 1993.
		
		\bibitem[SXY19]{SXY19}
		Z.~Shi, S.~Xu, and D.~Yan.
		\newblock Damping estimates for oscillatory integral operators with
		real-analytic phases and its applications.
		\newblock {\em Forum Math.}, 31(4):843--865, 2019.
		
		\bibitem[TX24]{TX24}
		Y.~Tan and S.~Xu.
		\newblock Some new decay estimates for (2+ 1)-dimensional degenerate
		oscillatory integral operators.
		\newblock {\em Arch. Math.}, 122(4):437--447, 2024.
		
		\bibitem[Var76]{Var76}
		A.~N. Varchenko.
		\newblock Newton polyhedra and estimation of oscillating integrals.
		\newblock {\em Funct. Anal. Appl.}, 10(3):175--196, 1976.
		
		\bibitem[Xia17]{Xia17}
		L.~Xiao.
		\newblock Endpoint estimates for one-dimensional oscillatory integral
		operators.
		\newblock {\em Adv. Math.}, 316:255--291, 2017.
		
		\bibitem[Xu23]{Xu23}
		S.~Xu.
		\newblock A sharp decay estimate for degenerate oscillatory integral operators
		using broad-narrow method.
		\newblock {\em J. Geom. Anal.}, 33(4):115, 2023.
		
		\bibitem[Yan04]{Yan04}
		C.~W. Yang.
		\newblock Sharp {$L^p$} estimates for some oscillatory integral operators in
		$\mathbb{R}^1$.
		\newblock {\em Illinois J. Math.}, 48(4):1093--1103, 2004.
		
	\end{thebibliography}

\end{document}